\documentstyle[12pt]{article}
\title{The Riemann hypothesis - an elementary  analytic approach based on 
complex Laplace transform}
\author{Andrzej M\c{a}drecki\thanks{Institute of Mathematics and
Computer Science, Wroc{\l}aw University of Technology, Wybrze\.ze
Wyspia\'nskiego 27, PL-50-370 Wroc{\l}aw, Poland}}
\newtheorem{th}{Theorem}
\newtheorem{re}{Remark}
\newtheorem{ex}{Example}
\newtheorem{lem}{Lemma}
\newtheorem{de}{Definition}

\newfont{\lll}{msbm10 scaled 1095}
\def\LA{\mbox{\lll \char65}}
\def\LC{\mbox{\lll \char67}}
\def\LN{\mbox{\lll \char78}}
\def\LR{\mbox{\lll \char82}}
\def\LT{\mbox{\lll \char84}}
\def\LZ{\mbox{\lll \char90}}
\def\LQ{\mbox{\lll \char81}}
\begin{document}
\maketitle

{\bf Abstract}. An elementary analytic proof of the famous Riemann
hypothesis is given. The main "accent" of the proof is a both using
of the 2-dimensional double real and complex Laplace integral
representations of the Green function $\mid z \mid^{-2}$.

\section{Introduction.}
Let $\LR$ and $\LC$ be the real and complex number fields, respectively.
In [ML] we gave a proof of the {\bf Riemann hypothesis}(Rh for short)
based on the following two-dimensional and two-sided {\bf Laplace
representation} of the complex {\bf Green function} $G_{2}(z)= \mid z
\mid^{-2}$ :
\begin{equation}
\mid z \mid^{-2}=\int_{\LC}e^{z \cdot c}dH_{2}^{*}(c)\;,z \in
\LC^{*}:=\LC-\{0\},
\end{equation}
where $z\cdot c = re(z)re(c)+im(z)im(c)$ is the usual Euclidean scalar
product in $\LC=\LR^{2}$, $H_{2}^{*}$ is the so called {\bf Hodge
measure}(see [ML] for the motivation of the name) and all in the sequel
$re(z), im(z)$ denotes the real and imaginary part of a complex number $z$.

To obtain the representation (1.1) we used in [ML] the machine of
stochastic analysis based on the Wiener measure $w_{\infty}$.

In [MR] we gave another two constructions of $H_{2}^{*}$ based on
l-adic Wiener measure $w_{l}$ and l-adic Gibbs measure $G_{l}$ (see [MR]
- for details).

In this paper, basing strongly on some {\bf Pitkannen's idea} (trial) of a
quantum mechanics (conformal) proof of Rh (see [Pi]), we will be able -
in a trivial way - establish the {\bf elementary complex Laplace representation}
$Rep(\LC)$ of $G_{2}$: let $\chi_{\LR_{+}}(x)$ be the characteristic function of
$\LR_{+}:=[0,+\infty)$ in $\LR$. It is so called {\bf Heaviside
function} - important in electricity (see e.g. [LZ, IV.3, Example 1]). Using 
its, we have the following {\bf trivial formula} 
\begin{equation}
\frac{1}{z}=\int_{\LR}e^{-zl}\chi_{\LR_{+}}(l)dl=\int_{0}^{\infty}e^{-zl}
dl={\cal L}_{\LC}(\chi_{\LR_{+}})(z)
\end{equation}
if $re(z)>0$. 

In fact the formula (1.2) is a little bit formal. The strict proof of
(1.2) requires the calculations of the real Laplace transforms of sine
and cosine. Thus we have : let $z=u+iv$. Then (see e.g. [LZ])
\begin{displaymath}
\int_{0}^{\infty}e^{-zl}dl= \int_{0}^{\infty}e^{-ul}cos
vl dl-i\int_{0}^{\infty}e^{-ul}sinvl dl =
\end{displaymath} 
\begin{displaymath}
={\cal L}(cosvl)(u)-i{\cal
L}(sinvl)(u)=\frac{u-iv}{u^{2}+v^{2}}=\frac{\tilde{z}}{\mid z \mid
^{2}}=\frac{1}{z}.
\end{displaymath}

According to the {\bf Pitkanen's patent} for $\mid z \mid^{2}$
and {\bf Fubini theorem} we get
\begin{equation}
\frac{1}{\mid z \mid^{2}}=\frac{1}{z\tilde{z}}=\int_{0}^{\infty}
e^{-zl_{1}}dl_{1} \int_{0}
^{\infty}e^{-\tilde{z}}l_{2}dl_{2}=
\end{equation}
\begin{displaymath}
=\int
\int_{\LR_{+}^{2}}e^{-[zl_{1}+\tilde{z}l_{2}]}dl_{1}dl_{2}=\int_{\LR_{+}^{2}}e
^{-<z,l>}d^{2}l,
\end{displaymath}
for each $z \in \LC$ with $re(z)>0$, where $c(x+iy)=\tilde{x+iy}:=x-iy$ is the
{\bf complex conjugation}.
(Let us explain what we mean here by the Pitkannen's patent : in [Pi],
he used the methods of the quantum conformal field theory. He tried to
realize the famous {\bf Hilbert-Polya conjecture} that complex zeroes
$\rho_{n}$ of $\zeta$ should be {\bf eigenvalues} of $1/2+iH$, where $H$
is a hypothetical {\bf hermitian operator} on a hypothetical Hilbert
space ${\cal H}$, i.e. $\rho_{n}=ev_{n}(iH+1/2)$ (with a discrete
spectrum), having the decomposition : $H=DD^{+}$ with such
non-hermitian operators $D$ and $D^{+}$ with the property that they
have zeroes $\rho_{n}$ and their conjugates as their complex
eigenvalues. Thus, we can write : $D(\psi_{n})=\rho_{n}\psi_{n}$ and
$D^{+}(\psi_{n})=\tilde{\rho_{n}}\psi_{n}$, for some eigenvectors
$\psi_{n} \in {\cal H}$).

Through this paper, we are using the bi-linear form $<.,.>:\LC \times
\LR \longrightarrow \LR$ given by (1.3), i.e. we put
\begin{equation}
<z,l>:=zl_{1}+c(z)l_{2},
\end{equation}
if $z \in \LC$ and $l=(l_{1},l_{2})\in \LR^{2}$.

The existence of the complex Laplace integral representation (1.3) (or the
double Fresnel integral (oscilatory) representation):
\begin{displaymath}
\mid z \mid^{-2}=\int
\int_{\LR_{+}^{2}}e^{-re(z)(l_{1}+l_{2})}cos(im(z)(l_{1}-l_{2}))dl_{1}dl
_{2}
\end{displaymath}
is a very happy circumstance (event) - it solves the so called - {\bf Zabczyk's
problem (question)} (Z.q. for short) : prove - in an elementary way -
that the function of the complex variable $\mid z \mid^{-2}$ is LHpd
(or {\bf positively definite} in the {\bf Toeplitz sense}(see [A]))
on $\LR_{+}^{2}$ - posed at a Seminar on Stochastic Analysis at IMPAN
in Warsaw - during the presentation by the author the so called Hodgian
proof of the Riemann hypothesis. (Let us mention that the same problem
was erected by {\bf M. Bo\.zejko} - by a private communication - many
years earlier. Thus, it is may be better to say on a {\bf
Bo\.zejko-Zabczyk question}(BZq - for short)).

Finally, we remark that similarly like in [MW] we proved (Rh) using the
methods of functional analysis and probability theory, in [MH] - the
methods of harmonic analysis and arithmetics of (Rh) - the methods of
that paper are concentrated (belong) to the Laplace transform theory
(important - for example in probabilistic and theory of partial
differential equations) and the classical elementary analysis. Thus it
is an elementary analysis of the Riemann hypothesis.

\section{The Fresnel integrals, Bernstein-Widder theorem and existence
of the Hodge measure $H_{2}$.}

In [AM] we introduced and considered some particularly useful and
simple property of Fourier integrals on $\LR$ : the {\bf non-vanishing}
of the so called {\bf Fresnel integrals} (cf.e.g. [S, IV.9., Example],
[EFI, p. 635], [PEGK]). For a very detailed and deep study of Fresnel
integrals see [AC]

Fresnel was probably the first, who considered the following
(fundamental in optics) {\bf Fresnel integrals} :
\begin{displaymath}
F_{2}(\nu):=\int_{0}^{\infty}sin \nu x^{2}dx=\int_{0}^{\infty}cos \nu
x^{2}dx=\frac{1}{2}\sqrt{\frac{\pi}{2\nu}},
\end{displaymath}
according to some problems concerning the light scattering theory.

In [S], {\bf L. Schwartz} considered the integrals :
\begin{displaymath}
F(r,\nu):=\int_{0}^{\infty}\frac{sin \nu x}{x^{r}}\;for\;r>0\;and\;\nu>0.
\end{displaymath}
In particular, the value $F(1,\nu)=\pi/2$ does not depend on $\nu$ and 
\begin{displaymath}
F(1/2,\nu)=\int_{0}^{\infty}\frac{cos \nu
x}{\sqrt{x}}dx=2\int_{0}^{\infty}cos \nu
x^{2}=2F_{2}(\nu)=\sqrt{\frac{\pi}{2\nu}}. 
\end{displaymath}

Therefore, in the sequel, it would be very convenient to introduce the
following integrals : each {\bf Positive Continuous Integrable and (strictly)
Decreasing} function $A : \LR_{+}\longrightarrow \LR_{+}$ - we call in the sequel
a {\bf PCID-amplitude} (see [AM]) and any $\nu \in \LR_{+}$ - the {\bf
frequency}. Then we can consider the integrals :
\begin{equation}
F_{s}(A,\nu):=\int_{0}^{\infty}A(x)sin \nu x
dx\;and\;F_{c}(A,\nu):=\int_{0}^{\infty}A(x)cos \nu x dx.
\end{equation}
According to the mentional above fundamental examples, it will be very
convenient to call such integrals - the {\bf Fresnel integrals} (
associated with an amplitude $A(x)$ and a frequency $\nu$).

({\bf M. Pluta} has communicated to the author that - in fact in optics -
important are rather Fresnel integrals being the functions of the upper
limit in the definite integral - by private communication. But for the purposes
of this paper we need only Fresnel integrals as numbers).

Obviously, $F_{s}(A,\nu)$ and $F_{c}(A,\nu)$ are thus nothing that the
{\bf real} and {\bf imaginary parts} of the {\bf Fourier (osscilatory)
integrals} (see e.g. [Ar]), i.e.
\begin{displaymath}
{\cal F}(A, \nu):=\int_{\LR} A(x)e^{i \nu x}dx,
\end{displaymath}
and , obviously $F_{c}(A,\nu)=re({\cal F}(A,\nu))$ and
$F_{s}(A,\nu)=im({\cal F}(A, \nu))$.

We have the following simple (but extremaly useful) analytic elementary
lemma :
\begin{lem}({\bf The Fresnel Lemma}.)

The Fresnel integrals of PCID-amplitudes are strongly positive, i.e.
for each $\nu>0$ and each PCID-amplitude $A(x)$ holds:
\begin{equation}
F_{s}(A,\nu)>0\;and\;F_{c}(A,\nu)=-\frac{1}{\nu}F_{s}(A^{\prime}, \nu).
\end{equation}
\end{lem}
{\bf Proof}. We have 
\begin{displaymath}
F_{s}(A,\nu)=\int_{0}^{\infty}A(x)sin(2\pi \nu
x)dx=\sum_{n=0}^{\infty}\int_{n/\nu}^{(n+1)/\nu}A(x)sin(2\pi \nu x)dx=
\end{displaymath}
\begin{displaymath}
=(\sum_{n=0}^{\infty}(-1)^{n}A(x_{n}))P(\nu),
\end{displaymath}
where $P(\nu)=\int_{0}^{1/\nu}sin 2\pi \nu x dx>0$ and a sequence
$\{x_{n}\}$ with $x_{n}\in [n/2\nu, (n+1)/\nu]$ is determined according
to the elementary mean value theorem (see [AM]).

Moreover, integrating by parts we have
\begin{displaymath}
F_{c}(A,\nu)=\int_{0}^{\infty}A(x)cos \nu x dx=(A(x)sin \nu x)/\nu
\mid_{0}^{\infty}-\frac{1}{\nu}\int_{0}^{\infty}A^{\prime}(x)sin \nu x
dx=
\end{displaymath}
\begin{displaymath}
=-\frac{1}{\nu}F_{s}(A^{\prime}, \nu).
\end{displaymath}

Since
\begin{equation}
\mid z \mid^{-2}=re(\int \int_{\LR_{+}^{2}}e^{-<z,l>}d^{2}l)=\int
\int_{\LR_{+}^{2}}
e^{-x(l_{1}+l_{2})}cos(y(l_{2}-l_{1}))dl_{1}dl_{2}=:F_{22}(z)\;,\;z=x+iy,
\end{equation}
thus $\mid z \mid^{-2}$ is nothing that the {\bf double Fresnel integral}.
Doing the change of variables in $F_{22}(z)$ (the Jacobi theorem),
according to the substitution : $l_{1}=u, l_{2}-l_{1}=v$, we  obtain
\begin{displaymath}
F_{22}(z)=\int
\int_{\LR_{+}^{2}}e^{-x(v+2u)}cos(yv)dudv=\int_{0}^{\infty}e^{-2xu}du
\int_{0}^{\infty}e^{-xv}cos(yv)dv=
\end{displaymath}
\begin{displaymath}
= F_{c}(exp^{-2re(z)})(0)F_{c}(exp^{-re(z)})(im(z))>0,
\end{displaymath}
according to the {\bf Fubini theorem}.

Let us recall that we have for the disposition the following two
{\bf bilinear forms} : the {\bf complex form}
\begin{displaymath}
<z,l> : \LC_{++}\times \LR_{+}^{2}\longrightarrow \LC,
\end{displaymath}
and the {\bf real form} (the Hilbert (Euclidean) scalar product) $z
\cdot l : \LC \times \LR_{+}^{2} \longrightarrow \LR$
\begin{displaymath}
z \cdot l :=re(z)l_{1}+im(z)l_{2}\;\;,\;\;l=(l_{1},l_{2}).
\end{displaymath}

For any positive $\sigma$-additive {\bf measure} $\mu$ on $(\LR_{+}^{2},
{\cal B}(\LR_{+}^{2}))$, where ${\cal B}(\LR_{+}^{2})$  is the $\sigma$-field
of all {\bf Borel subsets} of $\LR_{+}^{2}$, we can consider two {\bf
double Laplace transforms} : 

1. the {\bf complex Laplace transform} of a measure $\mu$,
\begin{equation}
{\cal L}_{2}^{+}(\LC)(\mu)(z):=\int \int_{\LR_{+}^{2}}
e^{-<z,l>}d\mu(l),
\end{equation}
and 

2. the {\bf real Laplace transform} of a measure $\mu$,
\begin{equation}
{\cal L}_{2}^{+}(\LR)(\mu)(z):=\int \int_{\LR_{+}^{2}}e^{-z \cdot
l}d\mu(l)=\hat{\mu}(z).
\end{equation}
\begin{re}
Let us recall (see [ML]) that we have in fact two (completely different)
harmonic analysis notions of the {\bf positive definity} (p.d. for short).

We say that a function $l : (S,+)\longrightarrow \LR_{+}$, defined on an 
{\bf abelian semigroup} $(S, +)$ is {\bf positive definite} in {\bf Laplace-Hankel 
sense} (LHpd in short - or simply we say on the {\bf semigroup positive-
definity}) iff for each real n-tuplet $(r_{1}, ... , r_{n})\in \LR^{n}$
and each semigroup n-tuplet $(s_{1}, ... , s_{n}) \in S^{n}$ holds
\begin{displaymath}
(LHpd)\;\;\sum_{i=1}^{n}\sum_{j=1}^{n}r_{i}r_{j}l(s_{i}+s_{j})\ge 0, 
\end{displaymath}
whereas, a function $f : (G, +)\longrightarrow \LC$, defined on an {\bf
abelian group} $(G, +)$ is called the {\bf positive definite} in the
{\bf Fourier-Hermite sense} (FHpd for short - or simply we say on the
{\bf group positive-definity}) iff for any n-tuples $g=(g_{1}, ... ,
g_{n})$ of elements of $G$ and $c=(c_{1}, ... , c_{n})$ of complex
numbers holds
\begin{displaymath}
(FHpd)\;\;\sum_{i=1}^{n}\sum_{j=1}^{n}c_{i}\tilde{c_{j}}f(g_{i}-g_{j})
\ge 0.
\end{displaymath}
\end{re}

The real Laplace transforms ${\cal L}_{2}^{+}(\LR)(\mu)=\hat{\mu}$ are
obviously LHpd (see also [ML] and [A]). But it is the very surprising
fact that the {\bf algebraically-ordered condition} LHpd gives the {\bf
characterization} of Laplace transforms on LCA-semigroups (see [ML] and
[A, Th. 5.5.4 and p.228]), according to the following deep
\begin{th}({\bf The Bernstein-Widder theorem for} $\LR^{n}$).

In order that the representation
\begin{equation}
l(x)\;=\;\int_{\LR^{n}}e^{x \cdot u}dh_{l}(u),
\end{equation}
be valid, where $h_{l}(u)$ is an {\bf unique non-decreasing function on
each factor variable}, it is necessary and sufficient that $l(x)$ is
{\bf continuous} LHpd on $\LR^{n}$.
\end{th}
\begin{re}
Let us mention the following beautiful {\bf Sierpi\'nski theorem} - at
this moment - if $f :(a,b) \longrightarrow \LR$  is such a function that
$-\infty< f(x) \le \infty (a <x<b)$, is convex and measurable, then
$f(x)$ is continuous.
\end{re}

\begin{th}({\bf On the existence of the Hodge measure} $H_{2}$).

For each $z \in \LC$ with $re(z)>0$ and $im(z)>0$, we have
\begin{equation}
Rep(\LR)\;\;\;\mid z \mid^{-2}=\int\int_{\LR_{+}^{2}}e^{-z\cdot l}dH_{2}(l)={\cal
L}_{2}^{+}(\LR)(H_{2})(z).
\end{equation}
\end{th}
{\bf Proof}. We show that the function $g_{2}(z)=\mid z \mid^{-2}$ (the
second Green function) is LHpd on $\LR_{+}^{2}$, i.e. for each
$z=(z_{1}, ... , z_{n})\in (\LR_{+}^{2})^{n}$ and $r=[r_{1}, ... ,
r_{n}]\in \LR^{n}$ holds :
\begin{equation}
\sum_{i=1}^{n}\sum_{j=1}^{n}r_{i}r_{j}\mid z_{i}+z_{j} \mid^{-2}\ge 0.
\end{equation}
In other words, the quadratic form :
\begin{displaymath}
r G_{2}(z)r^{T}:=[r_{1}, ..., r_{n}][\mid z_{i}+z_{j}\mid^{-2}]_{n
\times n}[r_{1}, ... , r_{n}]^{T},
\end{displaymath}
is LHpd in the {\bf Sylvester sense}(see e.g. [Ko]).

(Let us mention on the danger of the fact that $\mid 0 \mid^{-2}=+\infty$).
In this case we have deal with indefinite symbols : $\infty - \infty $.

Using the complex Laplace representation (1.3) we obtain (first we take elements
from the semigroups $S_{a}=\LR_{a}^{2}:=\{(x,y)\in \LR^{2} : x \ge a, y
\ge a\}, a>0$, i.e. $z \in S_{a}$)
\begin{equation}
\sum_{i=1}^{n}\sum_{j=1}^{n}r_{i}r_{j}\mid z_{i}+z_{j} \mid^{-2} \;=\;
\end{equation}
\begin{displaymath}
=\sum_{i=1}^{n}\sum_{j=1}^{n}r_{i}r_{j}\int
\int_{\LR_{+}^{2}}e^{-re(z_{i}+z_{j})(l_{1}+l_{2})} re(e^{i
im(z_{i}+z_{j})(l_{1}-l_{2})})dl_{1}dl_{2} \;=\;
\end{displaymath}
\begin{displaymath}
=\sum_{i=1}^{n}\sum_{j=1}^{n}r_{i}r_{j}\int
\int_{\LR_{+}^{2}}e^{-re(z_{i})(l_{1}+l_{2})}e^{-re(z_{j})(l_{1}+l_{2})}
[cos(im(z_{i})(l_{1}-l_{2}))cos(im(z_{j})(l_{1}-l_{2}))-
\end{displaymath}
\begin{displaymath}
-sin(im(z_{i})(l_{1})-l_{2}))sin(im(z_{j})(l_{1}-l_{2}))]dl_{1}dl_{2}=
\end{displaymath}
\begin{equation}
=\int\int_{\LR_{+}^{2}}(\sum_{i=1}^{n}\sum_{j=1}^{n}r_{i}cos(im(z_{i})
(l_{1}-l_{2})) r_{j}cos(im(z_{j})(l_{1}-l_{2}))
e^{-re(z_{i})(l_{1}+l_{2})}e^{-re(z_{j})(l_{1}+l_{2})})dl_{1}dl_{2} \;
\end{equation}
$+$
\begin{equation}
\int\int_{\LR_{+}^{2}}(\sum_{k=1}^{n}\sum_{j=1}^{n}[ir_{k}sim(im(z_{k})
(l_{1}-l_{2}))][-ir_{j}sin(im(\tilde{z}_{j})(l_{1}-l_{2}))] e^{re(z_{k}
(l_{1}+l_{2})}e^{-re(z_{j})(l_{1}+l_{2})})dl_{1}dl_{2}.
\end{equation}
Let us put :
\begin{equation}
r_{j}^{\prime}:=r_{j}cos(im(z_{j})(l_{1}-l_{2}))\in \LR,
\end{equation}
and
\begin{equation}
r_{j}^{\prime \prime}:=ir_{j}sin(im(z_{j})(l_{1}-l_{2}))\in \LC,
\end{equation}
$j=1, ... , n$. Then the first cosine double sum in (2.14) is equal to
\begin{equation}
\int\int_{\LR_{+}^{2}}(\sum_{j=1}^{n}r_{j}^{\prime}e^{-re(z_{j})(l_{1}+l
_{2})})^{2}dl_{1}dl_{2}\;\ge \;0,
\end{equation}
whereas, the second sine double sum in (2.15) is equal to
\begin{equation}
\int\int_{\LR_{+}^{2}}(\sum_{j=1}^{n}r_{j}^{\prime \prime}e^{-re(z_{j})
(l_{1}+l_{2})})c(\sum_{k=1}^{n}r_{k}^{\prime
\prime}e^{-re(z_{k})(l_{1}+l_{2})})dl_{1}dl_{2}\;=\;
\end{equation}
\begin{displaymath}
=\int\int_{\LR_{+}^{2}}\mid \sum_{j=1}^{n}r_{j}^{\prime
\prime}e^{-re(z_{j})(l_{1}+l_{2})}\mid^{2}dl_{1}dl_{2}\;\ge\;0.
\end{displaymath}

Applying the Bernstein-Widder theorem for $\mid z \mid^{-2}$ on
$\LR_{1/n}^{2}$ , we obtain the compatible sequence of {\bf finite measures}
$\{H_{2}^{n}\}$ with 
\begin{equation}
\mid z \mid^{-2}=\int \int_{\LR^{2}}e^{-z \cdot l}dH_{2}^{l}(l)\;,\;z \in
\LR_{1/n}^{2}.
\end{equation}
The inductive limit of $\{H_{2}^{n}\}$ gives obviously the required
measure $H_{2}$. 

Thus, we showed that $Rep(\LC)$ immediately implies $Rep(\LR)$, what is
very exciting, since for many years, it seemed that it is not possible,
and we always used the below Fernique-Haar measure systems.

Let us observe two "bad" properties of $H_{2}$ :

(i) $H_{2}(\LR^{2}) = +\infty$, i.e. $H_{2}$ is an infinite measure,

and

(ii) the {\bf support} $supp(H_{2})$ is $\LR^{2}_{+}$.

(Therefore, for many years- many maths people claimed that $H_{2}$
cannot exists! - by private communications) .

\section{$H_{2}$ is produced by any Fernique-Haar system $^{*}$.}

This section can be omitted under the first reading of the paper,
without any doubts for the understanding of this elementary analytic proof 
of the Riemann hypothesis.

According to the Bernstein-Widder theorem, we know that the measure
$H_{2}$ {\bf exists}, but we know nothing on the {\bf construction} (
or the inner structure) of the measure $H_{2}$, i.e. our proof of
existence of $H_{2}$ is not constructive - i.e. is "bad" - from the
point of view of {\bf Brouwer logic} and intuitionistic mathematics -
like the famous Cantor's proof of the existence of transcendental
numbers.

The below considered Fernique-Haar systems permit to {\bf construct}
$H_{2}$ and explain its inner structure.

Let $A$ be any {\bf measurable commutative algebra} with unit endowed with a
$\sigma$-{\bf field} ${\cal A}$ of subsets of $A$ and $I$ be a
{\bf measurable linear subspace} of $A$. Let $<. , .> : I \times A
\longrightarrow \LR$ be a $\LQ$-{\bf bilinear-measurable form} (w.r.t.
${\cal A}$) (see also [ML] and [MR] - for the different constructions
of Riemann-Weil cohomologies).

Moreover, we also assume that is given a {\bf way of immersion} of the
vector space $I$ into $A$ via the {\bf multiplication} on some non-zero
element $i_{0} \in I$, i.e. $M_{i_{0}}: I\longrightarrow A$, where
$M_{i_{0}}(i):=i\cdot i_{0}, i \in I$.

Moreover, we assume that there exists $I$-{\bf invariant probability
measure} $H : {\cal A} \longrightarrow [0,1]$ (in this paper -
similarly like in [MR] - we call such measures the {\bf weak (or residual) 
Haar measures}), since for each $i \in I$ and each $B \in {\cal A}$ we have :
\begin{displaymath}
(Inv)\;\;\;\;H(B+i)\;=\;H(B).
\end{displaymath}

Finally, we assume that there exists such a {\bf positive constant} $f
= f(A)>0$ that for each $m \in (0, f)$ the integral
\begin{displaymath}
(FC)\;\;\;\;F_{H}^{m}(A)\;:=\;\int_{A}e^{\frac{m<a^{2},1>}{2}}dH(a)
\end{displaymath}
is {\bf finite}. From (FC) we immediately obtain the existence of a family of
{\bf finite Feynmann-Kac measures} $F_{H}^{m}$, defined as :
\begin{equation}
F_{H}^{m}(B):= \int_{B}e^{m<a^{2},1>/2}dH(a)\;,\;B \in {\cal A}.
\end{equation}

In the sequel such systems (sixtets)
\begin{displaymath}
(A, {\cal A}, I, <. , .>, H, f(A) ) ,
\end{displaymath}
we call the {\bf Fernique-Haar systems}. In the proof of the Riemann
hypothesis such systems plays such a good role like the famous
{\bf Kolygavin-Euler systems} in the proof of the Fermat Last Theorem.

Let $\LA=(A, {\cal A},<.,.>, I, H, f(A))$ be any Fernique-Haar system.
Let $m \in (0, f(A))$ and let
\begin{displaymath}
dF_{m}(a)=e^{\frac{m<a,a>}{2}}dH(a),
\end{displaymath}
be the {\bf finite Feynmann-Kac measure}. Let us calculate the {\bf
Laplace transform} ${\cal L}_{I}(F_{m})$ of $F_{m}$ (w.r.t. the
$\LR$-linear-biform $m<.,.>$ and restricted to $I$):
\begin{displaymath}
{\cal
L}_{I}(F_{m})(i):=\int_{A}e^{m<ia,1>}dF_{m}(a)=\int_{A}e^{m<ia,1>}
e^{\frac{m<a^{2},1>}{2}}dH(a)=
\end{displaymath}
\begin{displaymath}
=e^{-\frac{m<i^{2},1>}{2}}\int_{A}e^{(m/2)(<i^{2},1>+2<ia,1>+<a^{2},1>)}d
H(a)=e^{-m<i,i>/2}\int_{A}e^{\frac{m<(i+a)^{2},1>}{2}}dH(a)=
\end{displaymath}
\begin{displaymath}
e^{-m<i^{2},1>/2}\int_{A}e^{m<a^{2},1>/2}dH(a)=e^{-\frac{m<i^{2},1>}{2}}F
_{m}(A),
\end{displaymath}
since, obviously $H$ is $I$-{\bf invariant} and the biform $<.,.>$ is 
$\LQ$-linear. (Let us mention at this place, that opposite to the von Neumann
-Weil theorem, can exists diferent (up to a constant) I-invariant measures if
$I \ne A$).

Let us fix any $i_{0}\in I \ne \{0\}$. Then according to the above
calculations, for any $u, v \in \LQ$ we have
\begin{equation}
{\cal L}_{I}(F_{m})(ui_{0}){\cal L}_{I}(F_{m})(vi_{0})
=e^{-m<(u^{2}+v^{2})i_{0}^{2},1>}F_{m}^{2}(A)=
\end{equation}
\begin{displaymath}
=\int_{A^{2}}e^{m(u<i_{0},a_{1}>+v<i_{0},a_{2}>)}d(F_{m}\otimes
F_{m})(a_{1}, a_{2}).
\end{displaymath}
Since - obviously - $\LQ$ is {\bf dense} in $\LR$, then the formula (3.22)
holds for all $(u,v) \in \LR^{2}$.

Integrating (3.22) with respect to the {\bf Lebesgue measure} $dt$ on
$\LR_{+}$ we get
\begin{equation}
\frac{F_{m}^{2}(A)}{m(u^{2}+v^{2})<i_{0}^{2},1>}
=F_{m}^{2}(A)\int_{\LR_{+}}e^{-[(u^{2}+v^{2})m<i_{0}^{2},1>]t}dt=
\end{equation}
\begin{displaymath}
=F_{m}^{2}(A)\int \int_{\LR_{+}\times
A^{2}}e^{[m\sqrt{t}(u<i_{0}a_{1},1>+v<i_{0}a_{2},1>)]}dt\otimes
d(F_{m}\otimes F_{m})(a_{1},a_{2}).
\end{displaymath}
Let us consider the measure space $(\LR_{+}\times A^{2}, dt \otimes
dF_{m}^{2})$ and the random variable 
\begin{displaymath}
X_{mi_{0}}(t,a_{1},a_{2})=m\sqrt{t}(<i_{0},a_{1}>, <i_{0},a_{2}>)\in
\LR^{2}.
\end{displaymath}
Let us put :
\begin{equation}
H_{2}^{*}:=\frac{m<i_{0}^{2},1>}{F_{m}^{2}(A)}X_{mi_{0}}^{*}(F_{m}^{2}
\otimes dt).
\end{equation}
Using the transport of measure theorem, we finally get : $z=(u,v)\in
\LC$
\begin{equation}
\mid z \mid^{-2}=\int\int_{\LR_{+}^{2}}e^{ux_{1}+vx_{2}}dX^{*}_{mi_{0}}
(F_{m}^{2}\otimes dt)(x_{1},x_{2})\;=
\end{equation}
\begin{displaymath}
=\;\int_{\LR_{+}^{2}} e^{-z \cdot l}dH_{2}^{*}(l),
\end{displaymath}
since $supp(H_{2})=\LR_{-}$.
According to the {\bf uniqueness} of a Laplace representation, we
finally get : $H_{2}^{*}= H_{2}$, i.e. $H_{2}^{*}$ is a required Hodge
measure and moreover is produced by our Fernique-Haar system.

\begin{re}
Our Fernique-Haar systems introduced and considered here are modeled on
the following {\bf Fernique-Girsanov system} considered in [ML] :
\begin{displaymath}
(C_{0},{\cal C}, \LR \cdot L_{a}, <,>_{L^{2}}, w_{\infty},
f(w_{\infty})),
\end{displaymath}
where $C_{0}$ is the Frechet space of all continuous functions $f
:\LR_{+}\longrightarrow \LR$, with $f(0)=0,  L_{a}$ is a peak function, the bilinear
form $<,>_{L^{2}}$ is some difference of Hilbert space $L^{2}$-scalar product
(with two diffent measures), $w_{\infty}$ is the {\bf standard Wiener
measure} and $f(w_{\infty})$ its Fernique's constant of $w_{\infty}$.

Obviously $w_{\infty}$ is not $\LR$-{\bf invariant}  but only
$\LR$-{\bf quasi-invariant}, since according to the {\bf Girsanov
theorem} we have
\begin{displaymath}
w_{\infty}(B+rL_{a}) =\int_{B}e^{-\frac{1}{2}r^{2}\int_{0}^{\infty}(
L_{a}^{\prime}(t))^{2}dt-r\int_{0}^{\infty}L_{a}^{\prime}(t)dx(t)}
dw_{\infty}(x), 
\end{displaymath}
where $r \in \LR$ and $B$ is a Borel set in $C_{0}$.

Let $C$ be the Frechet space of all real-valued continuous functions on
$\LR_{+}$. Then obviously we have the decomposition : $C=C_{0}\oplus
\LR$ and there exists a natural extension of $w_{\infty}$ to $C$, being
the distribution of the whole family of {\bf Brownian motions}
$B_{a}=(B_{a}(t): t \ge 0)$ starting from all points $a \in \LR$. Let
$B=(B_{a}(t): t \ge 0, a \in \LR)$. Let $W_{\infty}:=Law(B)$ be the
{\bf distribution (law)} of the stochastic field $B$. Unfortunately,
$W_{\infty}$ is not $\LR$-{\bf invariant}, since for each $a \in \LR$
the measure $W_{\infty}(\cdot + a)$ is {\bf singular} to $W_{\infty}$.
Realy, the distribution of a r.v. $B_{a}(0)$ is the Dirac delta point
measure $\delta_{a}$, i.e. $Law(B_{a}(0))=\delta_{a}$ and
$supp W_{\infty}(\cdot+a)\cap supp W_{\infty}=O$.

But, we can easy delete that disadventage, giving the following easy
$\LR$-invariant extension of $w_{\infty}$ to $C$ : since obviously
$C_{0}\oplus \LR \simeq C_{0}\times \LR$ then we can define the product
measure ${\cal W}_{\infty}$ by the formula :
\begin{displaymath}
{\cal W}_{\infty}:=w_{\infty}\otimes l,
\end{displaymath}
where $l$ is the Lebesgue measure on $(\LR,+)$. Then obviously ${\cal
W}_{\infty}$ is $\LR$-invariant but {\bf infinite} and do not satisfies
the Fernique condition (FC).

But it suffices to take the Frechet group $C(\LR_{+},\LT)$ of all
continuous functions on $\LR_{+}$ with values in the 1-dimensional
torus $\LT$ (instead of $\LR$) and a Haar probability measure $H_{\LT}$
on its, to obtain a proper measure $W_{\infty}:= w_{\infty}\otimes
H_{\LT}$.
\end{re}
\begin{ex}({\bf The l-adic Wiener-Fernique-Haar systems}.)

Let $C_{l}=C(\LZ_{l}, \LQ_{l})$ be the l-adic Banach space of all
$\LQ_{l}$-valued continuous functions defined on $\LZ_{l}$ (with the
sup-norm). Let ${\cal B}_{l}$ be its Borel $\sigma$-field. In [MR] we
constucted (combinatorically) the non-trivial $\LQ$-linear homomorphism
$I_{l}\in Hom_{\LQ}(\LQ_{l}^{+},
\LR^{+})=L_{\LQ}(\LQ_{l}^{+},\LR^{+})$. If $f, g \in C_{l}$, then we
define the biform $<.,.>_{l}$ by the formula :
\begin{displaymath}
<f,g>_{l}:=I_{l}(f(1)g(1)).
\end{displaymath}
Let $w_{l}$ be the standard l-adic Wiener measure  on ${\cal B}_{l}$(see
[MR]). Then $w_{l}$ has a compact $\LZ_{l}$-module support :
$K_{l}=supp(w_{l})$ and is $K_{l}$-invariant. Then the system :
\begin{displaymath}
(C_{l}, {\cal B}_{l}, <.,.>_{l}, K_{l}, w_{l}, +\infty),
\end{displaymath}
is a Fernique-Haar system. In particular we have got the l-adic
Feynmann-Kac measure $F_{l}$ given by the formula :
\begin{displaymath}
F_{l}(B)=\int_{B\cap K_{l}}e^{I_{l}(x^{2}(1))} dw_{l}(x)\;,\; B \in
{\cal B}_{l},
\end{displaymath}
(which do not depend on any parameter $m$). Since $w_{l}$ is a {\bf
Haar measure on} $K_{l}$(see [M\c{a}drecki, PhD. Thesis]) and
$I_{l}:\LQ_{l}\longrightarrow \LR$ is {\bf continuous}, then the
Fernique lemma is trivially satisfies in this case - and -
$f(C_{l})=f(w_{l})=+\infty$.

It is easy to see that the l-adic Laplace transform ${\cal L}(F_{l})$
has the form :
\begin{displaymath}
{\cal L}(F_{l})(x)=e^{-I_{l}(x^{2}(1))/2}F_{l}(C_{l}).
\end{displaymath}
\end{ex}
\begin{ex}({\bf The l-adic Gibbs-Fernique-Haar systems}).

Let $\LQ_{l}$ be the l-adic number fiel, $b_{l}$ its Borel
$\sigma$-algebra and $>.,<_{l}$ the bilinear form defined by the
formula :
\begin{displaymath}
>x,y<_{l}:=I_{l}(xy)\;,\;x, y \in \LQ_{l}.
\end{displaymath}
Finally, let $H_{l}$ be the Haar measure of $\LQ_{l}$ normalized by
$H_{l}(\LZ_{l})=1$. The the systems
\begin{displaymath}
(\LQ_{l}, b_{l}, >.,.<_{l}, l^{n}\LZ_{l}, H_{l}, +\infty),
\end{displaymath}
are the Fernique-Haar systems. We can define the l-adic Gibbs-F-K
-measures according to the formula :
\begin{displaymath}
f_{l}(B):=\lim_{n}\int_{B \in l^{n}\LZ_{l}}e^{I_{l}(x^{2})/2}dH_{l}(x),
\end{displaymath}
for {\bf bounded} set $B \subset l^{m}\LZ_{l}$ ( it is so called l-adic
termodynamical limit) and obviously there is the formula for suitable
l-adic (one-dimensional) Laplace transform of the form :
\begin{displaymath}
L(f_{l})(x)=e^{-I_{l}(x^{2})/2}f_{l}(l^{m}\LZ_{l}).
\end{displaymath}
\end{ex}
\begin{ex}{\bf The measurable cardinals-Pelc-Fernique-Haar systems}.

The notion of a measure was introduced to mathematics at the begining of
the twenty century in connection of some problems of real functions
theory. The measure theory leaded to some purely set theoretic problems.
One of such branches of the modern set theory is the theory of {\bf
measurable cardinals}.

In the modern set theory is also used the following extended notion of a
measure (see [K,Sect.10.6]).

If $X$ is a set and $f : X \longrightarrow \LR_{+}$ a positive
function, then by $\sum_{x \in X}f(x)$  is denoted the supremum of the
set of real numbers of the form $\sum_{x \in F}f(x)$, where $F$ is any
finite subset of $X$. In the case $X=\LN$ the sum $\sum_{x \in X}f(x)$
is equal to $\lim_{n}s_{n}$, where $s_{n}=\sum_{i<n}f(i)$.

For a set $X$, $\mid X \mid$ denotes the {\bf cardinality} of $X$ and
$P(X)$ the family of all subsets of $X$. Let $\kappa$ be any {\bf cardinal}.
We propose the following definition :
\begin{de}
Let $G$ be any {\bf discrete abelian group} and $I\triangleleft G$ a
non-zero {\bf subgroup} (Let us mention that any topological group can
be considered with its {\bf discrete topology} $G_{d}$ - see e.g.
[Hartman]). By a {\bf Pelc measure} $p$, we understand here any set
function which satisfies the following four conditions :

$(P_{1})$  $p : P(G)\longrightarrow [0,1], p(O)=0, p(G)=1$, i.e.
$p$ is a {\bf probability universal measure}.

$(P_{2})$  $p$ is $\kappa$-{\bf additive}, i.e. $p(\sum X)=\sum_{X \in {\cal X}}
p(X)$  for each family ${\cal X} \subset P(G)$
of commonly disjoint sets and such that $\mid {\cal X} \mid < \kappa$.
Here we consider only $\sigma$-additive measures (see [P]), i.e.
$\kappa = \omega = \mid \LN \mid$.

$(P_{3})$ $p(\{g\})=0$ for each $g \in G$, i.e. $p$ {\bf vanishes on
singletons}.

$(P_{4})$ $p$ is $I$-{\bf invariant}, i.e. for each $i \in I$ and any
subset $S$ of $G$ there is :
\begin{displaymath}
p(S+i)\;=\;p(S).
\end{displaymath}
\end{de}

Let us recal that a cardinal $\kappa$ is called a {\bf measurable
cardinal} (mc for short), if on some set of power $\kappa$ there exists
a universal probability measure (see [K] and [P, 0.1]).

Many classes of measures and different types of measurable cardinals -
e.g. real-valued and {\bf Ulam cardinals} were considered in [P], where
in particular, were established some beautiful theorems on the
existence of Pelc measures.

In this example we say shortly measurable cardinal (mc) like in [K] and
not like in [P], where is considered "the whole ZOO" of mc's.

As it is well-known, according to the famous {\bf Banach-Kuratowski
theorem}(see [BK]) : "if $2^{\omega}=\omega_{1}$ then $2^{\omega}$ is
not (mc)" and {\bf Ulam theorem}(see [P, Th.1.5]) - the {\bf general measure 
problem}, i.e. the problem : does there exist a $\sigma$-additive
$\sigma$-finite measure vanishing on singletons defined on all subsets
of the reals ? - cannot have a positive solution in usual mathematics.
In this setting the problem becomes purely set theoretic, i.e. depends
exclusively on the cardinality of the underlying set.

More exactly, assuming the {\bf Continuum hypothesis} (Ch for short),
such a measure cannot exist. But - as showed {\bf Cohn with Godel} - Ch
is {\bf independent} from the axioms of ZFC (see [K]). Nevertheless the
possibility of disproving the existence of (weakly){\bf inaccesible
cardinals} or real valued and Ulam or measurable cardinals is not
excluded. In spite of this danger the assumptions that such cardinals
exist are commonly used in modern set theory as {\bf additional
axioms}(see [P] and [K]), since :

$1^{0}$.  According to {\bf Solovay theorem}, the following statements
are equivalent :

$(S_{1})$ " (mc) exists" is consistent with ZFC

and

$(S_{2})$ The statement "$2^{\omega}$ is (mc)" is consistent with
ZFC,

as well as - according to the {\bf Kunen theorem}, the following are
equivalent :

$(K_{1})$ The statement "mc exists" is consistent with ZFC ,

$(K_{2})$ The statement " There exists a $\kappa$-complete
$\lambda$-saturated ideal on a cardinal $\kappa$" (where $\omega_{1}
\le \lambda \le \kappa^{+}$) is consistent with ZFC (see [P]).

Let us note that the positive answer to the original version of the
general measure problem (gmp), i.e. the statement "$2^{\omega}$ is mc " is in
a sense strictly stronger than the statement "there exists mc". {\bf
Levy} and {\bf Solovay} have proved that if the existence of (mc) is
consistent with set theory then the statement "there exists a
$mc+2^{\omega}=\omega_{1}$" is also consistent with ZFC and hence in
view of the above mentioned result of {\bf Banach} and {\bf Kuratowski}
the statement "there exists a (mc) $+2^{\omega}$ is not (mc)" is
consistent.

Finally a remark should be made about that set theoretic framework. All
the mentioned results are proved in usual set theory with the axiom of
choice(ZFC).

Similarly like (gmp) cannot be solved positively in the usual set
category, the {\bf invariant} and {\bf weak-invariant} gmp cannot be
positively solved in the usual (ZFC) category of {\bf groups},
according to two well-known {\bf deep negative results} (see [P]) of
{\bf Harazi$\tilde{s}$vili-Erdos-Mauldin} and a more general {\bf Ryll-
Nardzewski-Telgarsky'result}(see [RNT]) : if $I$ is any {\bf uncountable subgroup}
of a group $G$, there does not exists a {\bf universal} $I$-{\bf invariant}
$\sigma$-{\bf additive measure} on $G$, i.e. - in our terminology a
{\bf Pelc measure}.

Let $A$ be a commutative algebra with unit and endowed with such a $\LQ$-
bilinear form $<.,.>$ that it is {\bf exponentially square bounded} on
$A \times A$, i.e. 
\begin{displaymath}
sup_{a \in A}e^{<a,a>}<+\infty.
\end{displaymath}
Such algebras obviously exist - for example - we can take the compact ring
$\LZ_{l}$ for $A$ and the continuous biform $I_{l}(xy)$ for $<x,y>$.

Then obviously the {\bf Fernique Condition}(FC) is trivially satisfies
for any {\bf Pelc measure} $p : P(A)\longrightarrow [0,1]$, i.e. any
non-trivial universal probability and $A$-invariant measure $p$ ( it is
convenient to call such measures the {\bf Haar-Pelc measures}), since
\begin{equation}
F_{p}(A)=\int_{A}e^{\frac{<a^{2},1>}{2}}dp(a)\le sup_{a \in
A}e^{<a^{2},1>/2}<+\infty.
\end{equation}

Let a {\bf cardinal} $\alpha$ be the power of $A$, i.e. $\alpha =
\mid A \mid$. Assume that $\alpha$ is a {\bf measurable cardinal} (mc).
Then according to the {\bf Pelc theorem}[P, Th.2.5] - there exists a
{\bf Haar-Pelc measure} $p$ on A. Thus, we see that we have also to our
disposals the following {\bf measurable cardinals Pelc-Fernique
systems}:
\begin{equation}
(A, P(A), <.,.>, A=I, p , +\infty=:f(A)).
\end{equation}
The existence of systems (3.27) leads to the following - extremaly
surprising connection between the {\bf Continuum hypothesis}(Ch for
short) and the {\bf Riemann hypothesis}(Rh for short) given in the
following logical implication (let us call it the {\bf Measurable
Riemann hypothesis}(MRh for short)):

(MRh) {\bf The existence of measurable cardinals implies the Riemann
hypothesis} , i.e. $(mc)\Longrightarrow (Rh)$ (in short) (see also for
(mRh) in [ML] and the last remark of this paper).

Now, doing the above general simple quadratic calculus in this case of
measurable Pelc-Fernique-Haar systems we get
\begin{displaymath}
{\cal L}_{A}(p)(a)=e^{-<a^{2}/2,1>}F_{p}(A),
\end{displaymath}
i.e. we obtain the {\bf bounded Gaussian density} from the {\bf
unbounded} one.

Finally, let us remark, that if $G$ is a compact abelian group then
there exists the probability Haar measure on it, i.e. the unique
$G$-invariant probability measure $H_{G}$ defined on a very small
$\sigma$-algebra of {\bf Borel sets} ${\cal B}(G)$ of $G$. In
particular $p\mid {\cal B}(A)=H(A,+)$ is the probability Haar measure
if $A$ is compact.

If we take the {\bf Wiener measure} $w_{\infty}$ on the (non-locally compact)
group $G=C(\LR_{+})$ of all continuous real valued function on
$\LR_{+}$, then obviously $w_{\infty}(G)=1$ and obviously $p_{C(\LR_{+})}\ne
w_{\infty}$, since $w_{\infty}$ is not the Haar (or at least $\LR$-invariant). 
\end{ex}

\section{An elementary proof of the Riemann hypothesis.}

\begin{th}({\bf The Riemann hypothesis}).

Let $\zeta=\zeta(s)$ be the {\bf Riemann zeta function}. Then
\begin{displaymath}
(Rh)\;\;\;If\;\zeta(s)=0\;and\;im(s)\ne \;0\;then\;re(s)=\frac{1}{2}.
\end{displaymath}
\end{th}
{\bf Proof}. In this proof we will use the following tools from the
classical analysis :

1. The Riemann analytic continuation equation (Race for short).

2. The convergence of Dirichlet series.

3. The Newton-Leibnitz formula.

4. The change of summation and integration, i.e. the Fubini theorem for
series.

5. The $Rep(\LR)$ as well as $Rep(\LC)$ representations of the Green
function $\mid z \mid^{-2}$.
6. The positivity of Fresnel integrals.

Thus - in some sense - our proof is based on the shoulders of maths
giants.

($P_{1}$)({\bf Elementary analicyties and Rhfe}- {\bf from Race to
Rhfe}).

According to the {\bf classical Riemann analytic functional
continuation equation} (cf.e.g. [KV, L.2]) we have
\begin{equation}
i(s):=Im(\pi^{-s/2}\Gamma(\frac{s}{2})\zeta(s))=:Im(\zeta^{*}(s))=
\end{equation}
\begin{displaymath}
=Im(\frac{1}{s(s-1)})
+Im(\int_{1}^{\infty}(x^{\frac{s-2}{2}}+x^{-\frac{(s+1)}{2}})\theta(x)dx)
\end{displaymath}
, where all in the sequel by $G(x)=e^{-\pi x^{2}}$ we mean the {\bf
canonical Gaussian function} and $\theta(x)$ stands for the {\bf
canonical Jacobi theta}:
\begin{equation}
\theta(x):= \sum_{n=1}^{\infty}e^{-\pi
n^{2}x}=\sum_{n=1}^{\infty}G(n\sqrt{x}).
\end{equation}

Let us denote :
\begin{equation}
J(s):=\int_{1}^{\infty}(x^{\frac{(s-2)}{2}}+x^{-\frac{(s+1)}{2}})\theta(x)
dx,
\end{equation}
and
\begin{equation}
2J_{n}(x):=\int_{1}^{\infty}(x^{\frac{(s-2)}{2}}+x^{-\frac{(s+1)}{2}})G(n
\sqrt{x})dx.
\end{equation}
Let us observe that we can remove the fractions from the above
subintegral expression if we change variables according to the
substitution : $x^{2}=t, dx=2tdt$. Then
\begin{equation}
2J_{n}(s)=2\int_{1}^{\infty}(t^{s-1}+t^{-s})G(nt)dt.
\end{equation}

Now, let us observe that we have the following approximation of $\mid
J(s)\mid$ :
\begin{equation}
\sum_{n=1}^{\infty}\int_{1}^{\infty}(x^{u-1}+x^{-u})G(nx)dx \le
\mid\mid G \mid\mid_{(m,0)}(:=sup _{x \in \LR}\mid G(x)x^{m}\mid)\times
\end{equation}
\begin{displaymath}
\sum_{n=1}^{\infty}\int_{1}^{\infty}\frac{\mid x^{u-1}+x^{-u}
\mid}{n^{m}x^{m}} dx \le 2 \mid\mid G\mid\mid_{(m,0)} \zeta(m)max_{1\le
x \le \infty}(x^{u-1},x^{-u})\int_{1}^{\infty}\frac{dx}{x^{m}} <\infty,
\end{displaymath}
if $u=Re(s)\in I:=[0,1]$(in fact if $u\le 1$) and $m>1$.

Thus (4.33) shows that the {\bf iterated integral}
\begin{displaymath}
\sum_{n=1}^{\infty}\int_{1}^{\infty}\mid x^{u-1}- x^{-u} \mid G(nx) dx
\end{displaymath}
is {\bf finite} if $m>1$ and $u=Re(s)\le 1$.

Since the {\bf Lebesgue measure} $dt$ on $(\LR,+)$ and the {\bf calculating 
measure} $dc$ on $(\LZ, +)$ are $\sigma$-{\bf finite}, then according
to the {\bf Fubini theorem} for positive functions (cf.e.g. [B,
Commentaries after Th.18.2 and Exercise 18.3]), the functions
$f_{s}(x,n):= \mid x^{Re(s)-1}-x^{-Re(s)}\mid G(nx)$ are {\bf flat-
integrable} with respect to the product measure $dt\otimes dc$ on
$[1,+\infty)\times \LN$. Then according to the Fubini theorem for
arbitrary integrable functions (i.e. the {\bf Tonelli theorem}) we have
\begin{equation}
Im(J(s))-2\int_{1}^{\infty}(x^{Re(s)-1}-x^{-Re(s)})sin(Im(s)logx)\theta(x)
dx=
\end{equation}
\begin{displaymath}
=2\sum_{n=1}^{\infty}Im(J_{n}(s)).
\end{displaymath}

The {\bf detailed calculation} of $Im(J_{n}(s))$ for $s$ with $Re(s)\in
(0,1)$ is a "heart" of this "complex Laplace representation proof of
(Rh)" - in the spirit of the {\bf classical analytic number theory}. 

That calculations and the effect of them was quite surprising for us
many years ago!

Thus, for each $n \in \LN$ and $s=u+iv$ with $Re(s)\in (0,1)$ and
$v=Im(s)\in \LR^{*}$ we have 
\begin{equation}
Im(J_{n}(s))=\int_{1}^{\infty}(x^{u-1}-x^{-u})sin(vlogx)G(nx)dx=
\end{equation}
\begin{displaymath}
\int_{0}^{\infty}(e^{r(u-1)}-e^{-ru})sin(vr)G(ne^{r})e^{r}dr,
\end{displaymath}
if we apply the {\bf change variables formula}: $x=e^{r}, dx=e^{r}dr$.

But
\begin{equation}
\int_{0}^{\infty}(e^{ru}-e^{r(1-u)})sin(vr)e^{-\pi n^{2}e^{2r}}dr=
\end{equation}
\begin{displaymath}
=\lim_{N \rightarrow
\infty}\int_{0}^{N}(\sum_{j=0}^{\infty}\frac{(-1)^{j}(\pi
n^{2}e^{2r})^{j}}{j!})(e^{ru}-e^{r(1-u)})sin(vr)dr.
\end{displaymath}
Since the {\bf Taylor expansion} of $G$ is uniformly convergent on each
closed segment of $\LR$, then
\begin{equation}
\int_{0}^{N}(e^{ru}-e^{r(1-u)})sin(vr)(\sum_{j=0}^{\infty}\frac{(-\pi
n^{2} e^{2r})^{j}}{j!})dr=
\end{equation}
\begin{displaymath}
\sum_{j=0}^{\infty}\frac{(-\pi
n^{2})^{j}}{j!}\int_{0}^{N}(e^{r(2j+u)}-e^{(2j+1-u)})sin(vr)dr.
\end{displaymath}

But, it is an elementary fact that the following defined integrals are
easy calculable by the {\bf Newton-Leibnitz formula}
\begin{equation}
\int_{0}^{N}e^{wr}sin(vr)dr=\frac{e^{wr}(w
sin(vr)-vcos(vr))}{w^{2}+v^{2}}\mid_{0}^{N}=
\end{equation}
\begin{displaymath}
e^{Nw}\frac{(wsin(vN)-vcos(vN))}{w^{2}+v^{2}} +\frac{v}{w^{2}+v^{2}}.
\end{displaymath}

Since always (for a fixed $v=Im(s)\ne 0$) we can choose the sequence in
(4.37) and (4.38) of the form :
\begin{equation}
N:=\frac{2\pi L}{v}\;,\;L \in \LN,
\end{equation}
then the formula (4.38) will obtain a simpler form 
\begin{equation}
\int_{0}^{2\pi L/v}e^{wr}sin(vr)dr=\frac{v}{v^{2}+w^{2}}(1-e^{2\pi
Lw/v}),\;L\in \LN.
\end{equation}

$(P_{2})$({\bf The vanishing of the Poissonian part} $P_{n}(s)$
{\bf according to} $Rep(\LC)$).

Now, it will be convenient to do the following digression and notation
:  let $a$ and $b$ be arbitrary real numbers. By $P_{a}^{b}$ we denote
a {\bf generalized} and {\bf signed Poisson} random variable with the
parameters $a, b$. Thus
\begin{displaymath}
mes(P_{a}^{b}=k)=e^{-a}\cdot \frac{b^{k}}{k!},\;\;k \in \LN.
\end{displaymath}
If $a=b=\lambda>0$, then we obtain the {\bf Poisson distribution} with
a parameter $\lambda$.

Let $p_{a}^{b}(x)=mes(P_{a}^{b}\le x)$ be the generalized and signed
{\bf Poisson distribution function}. Then its {\bf Escher characteristic
function} $\hat{P}_{a}^{b}$ (being nothing that the {\bf real two-sided
Laplace transform with the unbounded exponential kernel}) has the form :
\begin{displaymath}
\hat{P}_{a}^{b}(u):=\int_{0}^{\infty}e^{ux}dp_{a}^{b}(x)=e^{-a}\sum_{k=0}
^{\infty}\frac{e^{uk}b^{k}}{k!}=e^{-a+be^{u}}, \;\;u \in \LR.
\end{displaymath}

Combining the identities from (4.34) to (4.40) we obtain
\begin{equation}
Im(J_{n}(s))=\lim_{L\rightarrow \infty}\sum_{j=0}^{\infty}\frac{(-\pi
n^{2})^{j}v}{j!} (e^{\frac{2\pi L(2j+1-u)}{v((2j+1-u)^{2})^{2}+v^{2}}}-
e^{\frac{2\pi L(2j+u)}{v(2j+u)^{2}}+v^{2}})+
\end{equation}
\begin{displaymath}
\sum_{j=0}^{\infty}\frac{(-\pi n^{2})^{j}}{j!}\cdot \frac{v(2u-1)(4j+1)}
{\mid (s+2j)(s-2j-1) \mid^{2}}=:
\end{displaymath}
\begin{displaymath}
=:P_{n}(s)\;+\;\zeta_{t}(s)Tr_{CG}^{n}(s).
\end{displaymath}
In particular, we see that the limit $Tr_{CG}^{n}(s)$ exists!

Additionally, we remark that
\begin{equation}
Im(\frac{1}{s(s-1)})=\frac{\zeta_{t}(s)}{\mid
s(s-1)\mid^{2}}=:\zeta_{t}(s)t_{0}(s),
\end{equation}
i.e. the {\bf zero-polar term} $\frac{1}{s(s-1)}$ of $\zeta(s)$ is
entered into the consideration of the below trace sequence 
$\{t_{n}(s)\}$.

Now, using the complex Laplace representation $Rep(\LC)$ of $\mid z \mid^{-2}$, we
show that for each $n \in \LN^{*}$, the {\bf Poissonian term} $P_{n}(s)$
{\bf vanishes}, i.e.
\begin{equation}
P_{n}(s)=0\;for\;n\in \LN^{*}\;and\;s\;with\;Re(s)\in I.
\end{equation}
To do that let us denote
\begin{equation}
P_{n}^{0}(z):=\lim_{L}\sum_{j=0}^{\infty}
\frac{e^{\frac{2\pi L(2j+Re(z))}{Im(z)}}}{\mid 2j+z \mid^{2}}
\frac{(-\pi n^{2})^{j}}{j!}=:
\end{equation}
\begin{displaymath}
=:\lim_{L}(\sum_{j=0}^{\infty}e^{a(z)L}\frac{(e^{b(z)L})^{j}}{\mid
2j+z\mid^{2}}\frac{(-\pi n^{2})^{j}}{j!}),
\end{displaymath}
where $z=s$ or $z=1-s, u=Re(s), v=Im(s)\ne 0$ and $j, L \in \LN$.
Finally, we have introduced also notations :
\begin{displaymath}
a(z):=\frac{2\pi Re(z)}{v}>0\;\;and\;\;b(z):=\frac{4\pi}{v}>0.
\end{displaymath}
According to (1.3) and (1.4) we have
\begin{equation}
\frac{1}{\mid 2j+z \mid^{2}}=\int_{\LR_{+}^{2}}e^{-<2j+z,l>}d^{2}l.
\end{equation}
Therefore, we have ($z=s$ or $z=1-s$)
\begin{equation}
P_{n}^{0}(z)=\lim_{L}(e^{a(z)L}(\sum_{j=0}^{\infty}\frac{(-\pi
n^{2})^{j}}{j!})(e^{b(z)L})^{j})\int_{\LR_{+}^{2}}e^{-<2j+z,l>}d^{2}l=
\end{equation}
\begin{displaymath}
\lim_{L}(e^{a(z)L}\int \int{\LR_{+}^{2}}(\sum_{j=0}^{\infty}\frac{(-\pi
n^{2}e^{b(z)L})}{j!}e^{-2(l_{1}+l_{2})})^{j})e^{-<z,(l_{1},l_{2})>}dl_{1}
dl_{2}                        =
\end{displaymath}
\begin{displaymath}
=\lim_{L}(e^{a(z)}\int \int_{\LR_{+}^{2}})e^{-\pi
n^{2}e^{(b(z)L-2(l_{1}+l_{2}))}}e^{-<z,l>}dl_{1}dl_{2}.
\end{displaymath}

The limit transition in the second line of the above three lines (4.46) is
become valid (under fixed $(n, L, z)$), since the assumptions of the
{\bf Tonelli theorem} for the {\bf space measure} $dc \otimes
dl_{1}\otimes dl_{2}$ are satisfies :
\begin{displaymath}
\sum_{j=0}^{\infty}\int \int_{\LR^{2}_{+}} \frac{(\pi
n^{2}e^{b(z)L})^{j}}{j!}e^{-2(j+Re(z))(l_{1}+l_{2})}
dl_{1}dl_{2}=\sum_{j=0}^{\infty}\frac{(\pi n^{2}
e^{b(z)L})^{j}}{j!}\cdot \frac{1}{4(j+Re(z))^{2}}\le
\end{displaymath}
\begin{displaymath}
\le (4 sup_{j}(j+Re(z))^{-2})e^{\pi n^{2}e^{b(z)L}}<+\infty.
\end{displaymath}
Finally, introducing the last notation : $c_{n}=\pi n^{2}$, we can
write down $P_{n}(z)$ as the limit :
\begin{equation}
P_{n}^{0}(z):=\lim_{L}P_{n}^{0}(L,z),
\end{equation}
where
\begin{equation}
P_{n}^{0}(L,z):=\int
\int_{\LR_{+}^{2}}e^{-c_{n}e^{b(z)L-2(l_{1}+l_{2})+a(z)}e^{-Re(z)(l_{1}+
l_{2})}}dl_{1}dl_{2}.
\end{equation}
But using a very coarse approximation : $e^{x}\ge \frac{x^{r}}{r!}, r
\in \LN^{*}$, we obtain (in short)
\begin{equation}
0\le e^{aL-c_{n}e^{(bL-2(l_{1}+l_{2}))}}\le
\frac{r!e^{aL}e^{2r(l_{1}+l_{2})}}{c_{n}^{r}e^{rbL}} \le
\end{equation}
\begin{displaymath}
\le  \frac{r!}{\pi^{r}n^{2r}}e^{(a-br)L}e^{2r(l_{1}+l_{2})}.
\end{displaymath}
Therefore
\begin{equation}
0\le \mid P_{n}^{0}(L,z) \mid \le
\frac{r!}{\pi^{r}n^{2r}}\int\int_{\LR_{+}^{2}}e^{\frac{2\pi}{Im(z)}(Re(z)
-2r)L}e^{2r(l_{1}+l_{2})}\mid e^{-<z,(l_{1},l_{2})>} \mid dl_{1}dl_{2}=
\end{equation}
\begin{displaymath}
=\frac{r!e^{\frac{2\pi}{Im(s)}(re(z)-2r)L}}{\pi^{r}n^{2r}}(\int\int_{\LR_{+}
^{2}}e^{-(re(z)-2r)(l_{1}+l_{2})}dl_{1}dl_{2}
=(re(z)-2r)^{-2})=:p_{L}(z,n,r).
\end{displaymath}
Under a fixed $z$ with $re(z)\in (0,1)$ and $n \in \LN^{*}$, taking
$r\ge 1$, we obtain
\begin{equation}
0\le \mid P_{n}^{0}(z) \mid  \le \lim_{L}p_{L}(z,n,r)=0.
\end{equation}

To finish the proof of this part of the proof, it sufficces to observe
that : $P_{n}(s)=P_{n}^{0}(s)-P_{n}^{0}(1-s)=0$.

$(P_{3})$({\bf There existence of a moment representation of the trace
sequence} $\{t_{j}(s)\}$).

Before we start to prove the subsequence part of the proof, we introduce the 
{\bf trace sequence} $t(s)=\{t_{j}(s)\}$ of $\zeta$ by the following
formula :
\begin{equation}
t_{j}(s):=\frac{4j+1}{\mid (s+2j)(2j+1-s) \mid^{2}} \;\;,\;j \in \LN.
\end{equation}

Now, we remark, that the complex Laplace transform representation $Rep(\LC)$ 
of $ \mid z \mid^{-2}$, is {\bf too weak} to obtain the subsequence part of the
elementary analytic proof of the Riemann hypothesis. We will base very
strongly on there existence of the {\bf Hodge measure} $H_{2}$ and
respectible {\bf real Laplace representation} $Rep(\LR)$. We prove that it induces the 
following : for each $s \in \LC$ with $re(s)\in (1/2,1)$ and $im(s)< 0$ holds :
\begin{equation}
t_{j}(s)\;=\;\int_{0}^{1} x^{j}dh_{s}(x)\;,\;j=0,1,2, ..., 
\end{equation}
i.e. $t_{j}(s)$ is the $j$-{\bf th moment} of the {\bf Hodge measure}
$h_{s}$.
Realy, let us write down the sequence $t(s)=\{t_{j}(s)\}$ of the form
\begin{equation}
t_{j}(s)=\frac{1}{\mid \frac{(0.5-s)}{4j+1}+0.5 \mid^{2}} \cdot
\frac{1}{4j+1}\cdot \frac{1}{\mid s+2j \mid^{2}}.
\end{equation}
In such a way, the decomposition of $t_{j}(s)$ is the main point of
this elementary proof of Rh.

The most difficult and fundamental in fact, is the possibility of a
real Laplace representation of the {\bf first factor} in the above
product, and therefore we first start from a consideration of that
factor. According to (2.20) , we have
\begin{equation}
\mid z \mid^{-2}
=\int\int_{\LR_{+}^{2}}e^{-z\cdot l}dH_{2}(l)\;,\; z \in \LR_{+}^{2}.
\end{equation}
Hence
\begin{equation}
\frac{1}{\mid (0.5-s)/(4j+1)+0.5
\mid^{2}}\;=\;\int\int_{\LR_{+}^{2}}e^{-((0.5-s)/(4j+1))\cdot
l}e^{-0.5\cdot l}dH_{2}(l)\;=
\end{equation}
\begin{displaymath}
\;=\;\int\int_{\LR_{+}^{2}}e^{\frac{-(0.5-re(s))\cdot l_{1}}{4j+1}}
e^{\frac{-im(s)\cdot l_{2}}{4j+1}}e^{-0.5\cdot l_{1}}dH_{2}(l_{1},l_{2}) \;=\;
\end{displaymath}
\begin{displaymath}
\int\int_{\LR_{+}^{2}}(\sum_{n=0}^{\infty}\frac{((re(s)-0.5)l_{1})^{n}}
{n!(4j+1)^{n}})
(\sum_{k=0}^{\infty}\frac{(-im(s)l_{2})^{k}}{k!(4j+1)^{k}}{k!(4j+1)^{k}})e
^{-0.5l_{1}}dH_{2}(l_{1},l_{2})\;=\;
\end{displaymath}
(Let us observe that all here , i.e. in this - the last but one part -
of the proof of the Riemann hypothesis concerning the existence of the
family of trace Hodge measures $\{h_{s}\}$, we assume that $re(s)\in
(1/2,1]$ and $im(s)<0$, since it guarantees the {\bf integrability} of the
subintegral function w.r.t. $H_{2}$).

But according to the {\bf Bernstein-Hausdorff theorem} (see [F.,
VII.9, Exercise 6 and XIII.4, Th.1 and XIII.1, Example(b)]), we have 
\begin{equation}
r^{-l}=\int_{0}^{\infty}e^{-ru}(\frac{u^{l-1}du}{(l-1)!}=:dB_{l}(u))\;,l
\ge 1, l \in \LN,
\end{equation}
, where $B_{l}$ is the "{\bf Bernstein measure}". (Let us observe that
in this place we have got a problem : the value $\Gamma(0)$ cannot be
represented by the integral $\int_{0}^{\infty}u^{l-1}e^{-ru}du$, for
$l=0$, since it last is divergent!).

Therefore, the integral in (4.56) we can write of the form :
\begin{equation}
\int_{\LR_{+}^{2}}(\sum_{k,n=0}^{\infty}\frac{((re(s)-0.5)l_{1})^{n}(-im
(s)l_{2})^{k}e^{l_{1}/2}}{n!(n-1)!k!(k-1)!})\times
\end{equation}
\begin{displaymath}
\times
\int_{0}^{\infty}e^{-(4j+1)u}u^{n-1}\int_{0}^{\infty}e^{-(4j+1)v}v^{k-1}d
v)dH_{2}(l_{1},l_{2}).
\end{displaymath}
Since the series $S(a):=\sum_{m=1}^{\infty}\frac{a^{m}}{m!(m-1)!}$ is
evidently absolutely convergent for all $a$, then we can change the
order of summazation and integration in the central integrals to obtain
the following form of the above integral :
\begin{equation}
\int_{\LR_{+}^{2}}dH_{2}(l_{1},l_{2})(\int\int_{\LR_{+}^{2}}(
\frac{e^{-4(u+v)})^{j}}{uv}\{((re(s)-1)ul_{1})\}
+\sum_{n=1}^{\infty}\frac{[(re(s)-0.5)ul_{1}]^{n}}{n!(n-1)!})e^{-l_{1}/2}
\times
\end{equation}
\begin{displaymath}
(\int_{0}^{\infty}(-im(s)vl_{2})+\sum_{k=1}^{\infty}
\frac{(-im(s)vl_{2})^{k}}{k!(k-1)!}dv)=
\end{displaymath}
\begin{displaymath}
=\int_{\LR_{+}^{2}}dH_{2}(l_{1},l_{2})\int\int_{\LR_{+}^{2}}(e^{-4(u+v)})^{j}
\frac{S((re(s)-0.5)ul_{1})S(-im(s)vl_{2})dudv}{e^{l_{1}/2+u+v}uv}.
\end{displaymath}
As usual - if we play a game with RH - to obtain an interesting
integral representation - we have to change an order of integration.
And here is a crucial moment of the using of the {\bf positivity} of the 
series $S(a)$. Then the above mentioned changing of integration is obvious - 
since the Fubini theorem always holds for the flat integrable positive function!.
Thus, we finally get for : $1/2<re(s)\le 1$ and $im(s)<0$):
\begin{equation}
\mid \frac{(0.5-s)}{4j+1}
+\frac{1}{2}\mid^{-2}=\int_{0}^{1}x^{j}dB_{s}^{1}(x),
\end{equation}
where $B_{s}^{1}(A):=\nu_{s}^{1}(\{(u,v)\in \LR_{+}^{2} :
e^{-4(u+v)}\in A\})$ (A is a Borel set in [0,1]) is the "distribution"
of the "random variable" $e^{-4(u+v)}$ on the measure space $(\LR_{+}^{2},
\nu_{s}^{1})$, where
\begin{equation}
d\nu_{s}^{1}(u,v):=(uv)^{-1}(\int\int_{\LR_{+}^{2}}e^{-l_{1}/2}S((re(s)-
0.5)ul_{1})S(-im(s)vl_{2})dH_{2}(l_{1},l_{2}))dudv.
\end{equation}

According to the Bernstein theorem (see (4.57)) applied in the case $l=1$,
there exists a positive measure $B^{2}$ with
\begin{equation}
\frac{1}{4j+1}\;=\;\int_{0}^{1}x^{j}dB^{2}(x).
\end{equation}
Finally, according to the existence of the real Laplace representation
of $\mid z \mid^{-2}$, there exists a positive measure $B_{s}^{3}$ with
the property :
\begin{equation}
\mid s+2j \mid^{-2} \;=\;\int_{0}^{1}x^{j}dB_{s}^{3}.
\end{equation}
Combining (4.60), (4.62) and (4.63) we finally get
\begin{equation}
t_{j}(s)=\int\int\int_{[0,1]^{3}}(xyz)^{j}dB_{s}^{1}(x)dB^{2}dB_{s}^{3}(x)
=\int_{0}^{1}x^{j}dh_{s}(x),
\end{equation}
where $h_{s}(A):=(B_{s}^{1}\otimes B^{2} \otimes
B_{s}^{3})(\{(x,y,z)\in [0,1]^{3} : xyz \in A\})$.

$(P_{4})$ ({\bf The finitness and strict positivity of zeta-traces}).

Let us introduce the following {\bf Cauchy-Gauss traces} of zeta (see
[ML]) :
\begin{equation}
Tr_{CG}(s)=\frac{1}{\mid s(s-1) \mid^{2}}\; +\;tr_{CG}(s),
\end{equation}
\begin{equation}
tr_{CG}(s):=\sum_{n=1}^{\infty}tr_{CG}^{n}(s),
\end{equation}
and
\begin{equation}
tr_{CG}^{n}(s):=\sum_{j=0}^{\infty}\frac{(-\pi n^{2})^{j}}{j!}
t_{j}(s).
\end{equation}
But according to (4.64) we have
\begin{equation}
tr_{CG}^{n}(s) = \sum_{j=0}^{\infty}\frac{(-\pi
n^{2})^{j}}{j!}\int_{0}^{1}x^{j}dh_{s}(x)=
\end{equation}
\begin{displaymath}
=\int_{0}^{1}(\sum_{j=0}^{\infty}\frac{(-\pi
n^{2}x)^{j}}{j!})dh_{s}(x)=\int_{0}^{1}e^{-\pi n^{2}x}dh_{s}(x),
\end{displaymath}
sine the {\bf Tonelli theorem} holds in this case. Indeed, let us
introduce the sequence of functions $\{f_{j}\}$ by the formula :
$f_{j}(x):=\frac{(\pi n^{2}x)^{j}}{j!}$. Then
\begin{equation}
0<\sum_{j=1}^{\infty}\int_{0}^{1}\mid f_{j}(x)\mid dh_{s}(x)\le
(\sum_{j=0}^{\infty}\frac{(\pi
n^{2})^{j}}{j!})sup_{j}(\int_{0}^{1}x^{j}dh_{s}(x)) =
\end{equation}
\begin{displaymath}
=e^{\pi n^{2}}sup_{j}h_{s}(j)<+\infty,
\end{displaymath}
since obviously the sequence $t(s)$ is {\bf bounded}. Thus, for each $n
\ge 1$
\begin{displaymath}
tr_{CG}^{n}(s)=\int_{0}^{1}e^{-\pi n^{2}x}dh_{s}(x)>0,
\end{displaymath}
and in a consequence we have : $Tr_{CG}(s)>0$.

On the other hand we have :
\begin{equation}
tr_{CG}^{n}(s)=:im(s)(tr_{CG}^{n0}(s)-tr_{CG}^{n0}(1-s)),
\end{equation}
where for $z=s$ or $z=1-s$
\begin{equation}
tr_{CG}^{n0}(z):= \sum_{j=0}^{\infty}\frac{(\pi n^{2})^{j}}{j!}\cdot
\frac{1}{\mid 2j+z \mid^{2}}=
\end{equation}
\begin{displaymath}
\sum_{j=0}^{\infty}\int_{\LR_{+}^{2}}\frac{(\pi
n^{2})^{j}}{j!}(e^{-2jl_{1}})^{j}e^{-z \cdot
l}dH_{2}(l_{1},l_{2})=\int_{\LR_{+}^{2}}e^{-\pi n^{2}e^{-2l_{1}}}e^{-z
\cdot l}dH_{2}(l).
\end{displaymath}
Therefore, applying the ineqality $e^{-x}\le \frac{d!}{x^{d}}$ for
arbitrary $d \in \LN^{*}, x>0$, we obtain
\begin{equation}
tr_{CG}^{0}:=\sum_{n=1}^{\infty}tr_{CG}^{n0} =
\end{equation}
\begin{displaymath}
=\sum_{n=1}^{\infty}\int_{\LR_{+}^{2}}e^{-\pi n^{2}e^{-2l_{1}}}e^{-z
\cdot l}dH_{2}(l)\le \sum_{n=1}^{\infty}\frac{d!}{(\pi
n^{2})^{d}}\int_{\LR_{+}^{2}}e^{-(z-2d)\cdot l}dH_{2}(l)=
\end{displaymath}
\begin{displaymath}
=\frac{d!}{\pi^{d}}\frac{\zeta(2)}{\mid z-2d \mid^{2}}<+\infty,
\end{displaymath}
if only $d \ge 1$ and $re(z) \in [0,1]$.

$(P_{5})$.{\bf The Riemann hypothesis functional equation}(Rhfe for
short).

Combining all calculations from the points $(P_{1}-P_{4})$- we obtain
the following {\bf Riemann hypothesis functional equation}(Rhfe in
short) for $re(s)\in [1/2, 1]$ and $im(s)<0$ : let
\begin{displaymath}
\zeta^{*}(s):=\pi^{-s/2}\Gamma(\frac{s}{2})\zeta(s),
\end{displaymath}
be the Gamma evolueted zeta. Then
\begin{equation}
im(\zeta^{*}(s))\;=\;im(s)(2re(s)-1)Tr_{CG}(s),
\end{equation}
with $Tr_{CG}(s)>0$, which obviously immediately implies the famous
{\bf Riemann hypothesis}.
\begin{re}
Let us observe that $\zeta(s)$ for $re(s)>1$ is nothing that the {\bf
complex Laplace transform} of the {\bf zeta measure} ${\cal Z}$, i.e.
\begin{displaymath}
\zeta(s)=\int_{0}^{\infty}e^{-sx}d{\cal
Z}(x)=\sum_{n=1}^{\infty}(e^{logn})^{-s},
\end{displaymath}
where
\begin{displaymath}
{\cal Z}\;:=\;\sum_{n=1}^{\infty}\delta_{logn},
\end{displaymath}
and obviously $\delta_{a}$ is the {\bf Dirac delta measure} at $a$.

In particular the {\bf Riemann hypothesis} (Rh for short) is {\bf
equivalent} to the following {\bf measure Riemann hypothesis}
(mRh) : there is a surprising, deep and enexpected relation between
the zeta measure ${\cal Z}$ with ${\cal L}({\cal Z})(s)=\zeta(s)$  and
the {\bf Lebesgue measure} $d^{2}l$ with ${\cal L}_{\LC}(d^{2}l)(z)= \mid z
\mid^{2}$ as well as the {\bf Hodge measure} $H_{2}$ with ${\cal
L}_{\LR}(H_{2})(s)=\mid s \mid^{-2}$ , i.e. formally,
\begin{displaymath}
{\cal Z} \longleftrightarrow (d^{2}l,dH_{2}).
\end{displaymath}
That connection is expressed explicite in the following {\bf
Cauchy-Gauss-Fresnel-Riemann hypothesis equation} (${\cal Z}d^{2}l$) :
\begin{displaymath}
({\cal
Z}d^{2}lH_{2})\;\;\;\;im(\pi^{-s/2}\Gamma(\frac{s}{2})\int_{\LR_{+}}e^{-sx}d{\cal
Z}(x))\;=\;
\end{displaymath}
\begin{displaymath}
\;=\;im(s)(2re(s)-1)\sum_{n=0}^{\infty}\sum_{j=0}^{\infty}\frac{(-\pi
n^{2})^{j}(4j+1)}{j!} \times
\end{displaymath}
\begin{displaymath}
\times(\int\int_{\LR_{+}^{4}}e^{-(<(2j+s),l_{1}>+<(2j+1-s,
l_{2})>}d^{2}l_{1}\otimes d^{2}l_{2}=\int\int_{\LR_{+}^{4}}
e^{-[(2j+s)\cdot l_{1}+(2j+1-s)\cdot l_{2}]}dH_{2}(l_{1})\otimes
dH_{2}(l_{2})\;,\;re(s)>1.
\end{displaymath}
Let us also remark that the sequence $\{\frac{(4j+1)}{j!}\}$ is quite
special! It appears as the sequence of coefficients in the Taylor
expansion of the {\bf canonical second Hermite function}.

Moreover, there exists also the connection
\begin{displaymath}
{\cal Z}\longleftrightarrow \{h_{s}\},
\end{displaymath}
given by the functional equation :
\begin{displaymath}
({\cal Z}h_{s})\;\;\;im(\pi^{-\frac{s}{2}}\Gamma(\frac{s}{2})\zeta(s))=
im(s)(2re(s)-1)\sum_{n=0}^{\infty}\sum_{j=0}^{\infty}\frac{(-\pi
n^{2})^{j}}{j!}\int_{0}^{1}x^{j}dh_{s}(x)=
\end{displaymath}
\begin{displaymath}
=im(s)(2re(s)-1)\sum_{n=0}^{\infty}\int_{0}^{1}e^{-\pi n^{2}x}dh_{s}=
im(s)(2re(s)-1)\int_{0}^{1}(\theta(x)+1)dh_{s}(x)=
\end{displaymath}
\begin{displaymath}
=:\zeta_{t}(s)\int_{0}^{1}(\theta(x)+1)dh_{s}(x)
\end{displaymath}
, where $\zeta_{t}(s):=im(s)(2re(s)-1)$ is the so called {\bf trivial
zeta}.
\end{re}

e-mails of the author : madrecki@o2.pl and madrecki@im.pwr.wroc.pl

\begin{thebibliography}{23}
\bibitem{1}[A] Akhiezer N.I., {\it The classical moment problem \/},
University Mathematical Monographs, Oliver and Boyd, Edinburgh and
London, 1965.
\bibitem{2}[AC] Albeverio S. and Cebulla C., {\it Muntz formula and
zero free regions for the Riemann zeta function \/},(preprint 2005),
p.1-27.
\bibitem{3}[AM] Albeverio S. and M\c{a}drecki A., {\it Probabilistic
proof of Riemann hypothesis\/}, (preprint 2006), p.1-41
\bibitem{4}[Ar] Arnold W.I.,Varchenko A.N. and Hussein-Zade S.M., {\it
 Singularities of differentials maps (Monodromy and asymptotic
integrals)\/}(in Russian), Moskou, Nauka, 1984.
\bibitem{5}[B] Billingsley P., {\it Probability and Measure \/}, John
Wiley and Sons, NY-Ch-B-T,1979.
\bibitem{5a}[BK] Banach S. and Kuratowski K., {\it Sur une
generalisation du probleme de la measure \/}, Fund. Math. 14(1929),
p.127-131.
\bibitem{6}[EFI] {\it Encyklopedia Fizyki I \/}, PWN , Warszawa 1972.
\bibitem{7}[F] Feller W., {\it An introduction to probability theory
and its applications II \/}, John Wiley and Sons Inc., NY 1966. 
\bibitem{8}[H] Hartman S., {\it Wst\c{e}p do Analizy Harmonicznej \/},
BM 33, PWN , Warszawa 1969.
\bibitem{9}[K] Kuratowski K. and Mostowski A., {\it Teoria Mnogo\'sci
\/}(in Polish) MM 27, PWN, Warszawa 1978.
\bibitem{10}[KV] Karatsuba A.A. and Voronin S.M., {\it The Riemann
Zeta-function\/}, de Gruyter Expositions in Mathematics 5, Berlin-NY,
Walter de Gruyter 1992.
\bibitem{11}[Ko] Komorowski J., {\it Od liczb zespolonych do tensor\'w,
spinor\'ow, algebr Liego i kwadryk  \/}(in Polish), PWN, Warszawa 1968.
\bibitem{12}[LZ] Leksi\'nski W. and \.Zakowski W., {\it Matematyka IV
\/}(in Polish), PWN, Warszawa 1976.
\bibitem{13}[MH] M\c{a}drecki A., {\it Hermitian proof of the Riemann
hypothesis \/}, (preprint 2005, reviewed in Mathematische Zeitschrift(MZ)
, p.1-32).
\bibitem{14}[ML] M\c{a}drecki A., {\it The Riemann hypothesis and some
stochastic Laplace reprezentation \/}(preprint 2005), reviewed in
Fundamenta Mathematica(FM), p.1-49.
\bibitem{15}[MR] M\c{a}drecki A., {\it The Riemann hypothesis for
Ramanujan zetas \/}, (preprint 2006), refereed in Manuscripta
Mathematica (MM), p.1-45.
\bibitem{16}[MW] M\c{a}drecki A., {\it A short Wiener measure proof of
the Riemann hypothesis \/}, accepted for its publication in the Journal
of Stochastic Analysis and Applications(JSAA), 30 November(2006), p.1-12.
\bibitem{17}[PEGK] Pluta M., Every A.G., Grill W., and Kim T.J. {\it
Fourier inversion of acoustic wave fields in anisotropic solids \/},
Physical Review B 67, 094117(2003), p.1-9.
\bibitem{18}[Pi] Pitkanen M. , {\it Proof of Riemann hypothesis \/}, 
(preprint 2001) : arXiv : math .GM/0102031, 5Feb(2001).
\bibitem{19}[P] Pelc A., {\it Invariant measures and ideals on discrete
groups \/}, Dissertationes Mathematicae CCLV, IM PAN, PWN, Warszawa
1986, p.1-45.
\bibitem{20}[RNT] Ryll-Nardzewski Cz. and Telgarsky R., {\it The
nonexistence of universal invariant measures \/}, Proc. Amer. Math.
Soc.(PAMS)69(1978), p.240-242.
\bibitem{21}[S] Schwartz L., {\it Analyse Mathematique I \/}(in French),
Cours profesee a l'Ecole Polytechnique, Hermann, Paris 1967.
\bibitem{22}[SZ] Saks S. and Zygmund A., {\it Funkcje Analityczne
\/}(in Polish), MM 10, PWN, Warszawa 1959.
\end{thebibliography}
\end{document}